\documentclass[12pt,leqno]{article}
\usepackage{amsfonts}
\usepackage{amsmath, amsthm, amsfonts, amssymb, color}
\usepackage{mathrsfs}
\usepackage{color}
\usepackage[notcite,notref]{showkeys}

\theoremstyle{definition}

\newcommand{\scr}[1]{\mathscr #1}
\definecolor{wco}{rgb}{0.5,0.2,0.3}

\numberwithin{equation}{section} \theoremstyle{remark}

\newcommand{\ua}{\uparrow}

\title{{\bf    Mild Solutions and Harnack Inequality for Functional SPDEs with Dini Drift  }\footnote{Supported in
 part by  NNSFC(11131003, 11431014), the 985 project and the Laboratory of Mathematical and  Complex Systems.} }
\author{
{\bf      Xing Huang $^{a)}$, Shao-Qin Zhang $^{b)}$}\\
\footnotesize{ a)School of Mathematical Sciences,
Beijing Normal
University, Beijing 100875, China,}\\
\footnotesize{  XingHuang@mail.bnu.edu.cn}\\
 \footnotesize{ b)School of Statistics and Mathematics, Central University of Finance and Economics, Beijing 100081, China, }\\
\footnotesize{  zhangsq@cufe.edu.cn }}
\begin{document}
\allowdisplaybreaks
\def\R{\mathbb R}  \def\ff{\frac} \def\ss{\sqrt} \def\B{\mathbf
B}
\def\N{\mathbb N} \def\kk{\kappa} \def\m{{\bf m}}
\def\ee{\varepsilon}\def\ddd{D^*}
\def\dd{\delta} \def\DD{\Delta} \def\vv{\varepsilon} \def\rr{\rho}
\def\<{\langle} \def\>{\rangle} \def\GG{\Gamma} \def\gg{\gamma}
  \def\nn{\nabla} \def\pp{\partial} \def\E{\mathbb E}
\def\d{\text{\rm{d}}} \def\bb{\beta} \def\aa{\alpha} \def\D{\scr D}
  \def\si{\sigma} \def\ess{\text{\rm{ess}}}
\def\beg{\begin} \def\beq{\begin{equation}}  \def\F{\scr F}
\def\Ric{\text{\rm{Ric}}} \def\Hess{\text{\rm{Hess}}}
\def\e{\text{\rm{e}}} \def\ua{\underline a} \def\OO{\Omega}  \def\oo{\omega}
 \def\tt{\tilde} \def\Ric{\text{\rm{Ric}}}
\def\cut{\text{\rm{cut}}} \def\P{\mathbb P} \def\ifn{I_n(f^{\bigotimes n})}
\def\C{\scr C}      \def\aaa{\mathbf{r}}     \def\r{r}
\def\gap{\text{\rm{gap}}} \def\prr{\pi_{{\bf m},\varrho}}  \def\r{\mathbf r}
\def\Z{\mathbb Z} \def\vrr{\varrho} \def\ll{\lambda}
\def\L{\scr L}\def\Tt{\tt} \def\TT{\tt}\def\II{\mathbb I}
\def\i{{\rm in}}\def\Sect{{\rm Sect}}
\def\M{\scr M}\def\Q{\mathbb Q} \def\texto{\text{o}} \def\LL{\Lambda}
\def\Rank{{\rm Rank}} \def\B{\scr B} \def\i{{\rm i}} \def\HR{\hat{\R}^d}
\def\to{\rightarrow}\def\l{\ell}\def\iint{\int}
\def\EE{\scr E}\def\no{\nonumber}
\def\A{\scr A}\def\V{\mathbb V}\def\osc{{\rm osc}}\def\H{\scr H}
\def\BB{\scr B}\def\Ent{{\rm Ent}}

\maketitle

\begin{abstract} The existence and uniqueness of the mild solution for a class of functional SPDEs with  multiplicative noise and a locally Dini continuous drift are proved. In addition, under a reasonable condition the solution is non-explosive. Moreover, Harnack inequalities are derived for the associated semigroup under certain global conditions, which is new even in the case without delay.

\end{abstract} \noindent
 \vskip 2cm

\section{Introduction}
Recently, using Zvonkin type transformation and gradient estimate, Wang \cite{W} has proved the existence and uniqueness of the mild solution for a class of SPDEs with multiplicative noise and a locally Dini continuous drift. Following this, Wang and Huang \cite{WH} extend the results to a class of Functional SPDEs, where the drift without delay is assumed as Dini's continuity, and the delay drift is Lipschitzian under $\|\cdot\|_{\C_{\nu}}$, see the details in \cite{WH} . In this paper, we try to replace $L^2(\nu)$ norm in \cite{WH} with uniform norm (finite delay) or weighted uniform norm (infinite delay). Due to the technique reason, for example, the Fubini Theorem is unavailable in the present case, we need a stronger condition on the singular drift than that in \cite{W}, see {\bf{(a3)}} in the following.

Let $\left(\mathbb{H},\langle,\rangle,|\cdot|\right)$ and $\left(\mathbb{\bar{H}},\langle,\rangle_{\mathbb{\bar{H}}},|\cdot|_{\mathbb{\bar{H}}}\right)$ be two separable Hilbert spaces.
For any $r\in(0,\infty]$, let $\C= C((-\infty,0]\cap [-r,0];\mathbb{H})$. For all $\xi\in\C$, define
\begin{equation*}
\|\xi\|_\infty =\sup_{s\in(-\infty,0]\cap [-r,0]} (\e^{-s}1_{r=\infty}+1_{r<\infty})|\xi(s)|.
\end{equation*}
For any $f\in C((-\infty,\infty)\cap [-r,\infty);\mathbb{H})$, $t\geq 0$, let $f_{t}(s)=f(t+s), s\in (-\infty,0]\cap [-r,0]$, then $f_{t} \in \C$.

Let $W=(W(t))_{t\geq 0}$ be a cylindrical Brownian motion on $\mathbb{\bar{H}}$ with respect to a complete filtration  probability space $(\OO, \F, \{\F_{t}\}_{t\ge 0}, \P)$. More precisely, $W(\cdot)=\sum_{n=1}^{\infty}{\bar{W}^{n}(\cdot)\bar{e}_{n}}$ for a sequence of independent one dimensional Brownian motions $\left\{\bar{W}^{n}(\cdot)\right\}_{n\geq 1}$ with respect to $(\OO, \F,
\{\F_{t}\}_{t\ge 0}, \P)$, where $\{\bar{e}_{n}\}_{n\geq 1}$ is an orthonormal basis on $\mathbb{\bar{H}}$.

Consider the following functional SPDE on $\mathbb{H}$:
\beq\label{1.1} \d X(t)= A X(t)\d t+b(t,X(t))\d t+B(t,X_{t})\d t+Q(t,X(t))\d W(t),\ \ X_0=\xi\in\C,
\end{equation}
where $(A,\D(A))$ is a negative definite self-adjoint operator on $\mathbb{H}$, $B: [0,\infty)\times \C\to \mathbb{H}$, $b: [0,\infty)\times \mathbb{H}\to \mathbb{H}$ are measurable and locally bounded (i.e. bounded on bounded sets), and $Q: [0,\infty)\times \mathbb{H}\to \L\left(\mathbb{\bar{H}}; \mathbb{H}\right)$ is measurable, where $\L\left(\mathbb{\bar{H}}; \mathbb{H}\right)$ is the space of bounded linear operators from $\mathbb{\bar{H}}$ to $\mathbb{H}$.

Let $\|\cdot\|$ and $\|\cdot\|_{\mathrm{HS}}$ denote the operator norm and the Hilbert-Schmidt norm respectively, and let $\L_{\mathrm{HS}}\left(\mathbb{\bar{H}}; \mathbb{H}\right)$ be the space of all Hilbert-Schmidt operators from $\mathbb{\bar{H}}$ to $\mathbb{H}$. Let $A, B$ and $Q$ satisfy the following two assumptions:
\beg{enumerate}
\item[{\bf (a1)}] $(-A)^{\varepsilon-1}$ is of trace class for some $\varepsilon \in(0,1)$; i.e. $\sum_{n=1}^{\infty}{\lambda_{n}^{\varepsilon-1}}<\infty$ for $0< \lambda_{1}\leq \lambda_{2}\leq\cdots\cdots$ being all eigenvalues of $-A$ counting multiplicities.

\item[{\bf (a2)}] $B\in C([0,\infty)\times \C; \mathbb{H})$, $Q\in C([0,\infty)\times \mathbb{H};\L(\mathbb{\bar{H}}; \mathbb{H}))$ such that for every $t\geq 0$, $Q(t,\cdot)\in C^{2}(\mathbb{H};\L(\mathbb{\bar{H}};\mathbb{H}))$, and $(Q Q^{\ast})(t,x)$ is invertible for all $(t,x)\in[0,\infty)\times \mathbb{H}$. Moreover,
\beg{equation*}\beg{split}
\|[\nabla B(t,\cdot)](\xi)\|+\sum_{j=0}^{2}\left\|[\nabla^{j} Q(t,\cdot)](\xi(0))\right\|+\left\|(Q Q^{\ast})^{-1}(t,\xi(0))\right\|
\end{split}\end{equation*}
\end{enumerate}
is locally bounded in $(t,\xi)\in[0,\infty)\times \C$, where $\|[\nabla B(t,\cdot)](\xi)\|$ stands for the local Lipschitz constant of $B(t,\cdot)$ at $\xi$.

Next, to describe the sigularity of $b$, we introduce
\beg{equation*}\beg{split}
\D= \Big\{\phi: [0,\infty)\to [0,\infty) \text{ is increasing}, \phi^{2} \text{ is concave}, \int_0^1{\frac{\phi(s)}{s}\d s}<\infty\Big\}
\end{split}\end{equation*}
and
\beg{equation*}\beg{split}
\A = \left\{a\in\B((0,\infty);(0,\infty)), \int_{0}^{1}\sup_{i\geq 1}\frac{\lambda_i\e^{-\lambda_is}}{a(\lambda_i)}\d s<\infty \right\}.
\end{split}\end{equation*}
\paragraph{Remark 1.1} The condition $\int_0^1{\frac{\phi(s)}{s}\d s}<\infty$ is well known as Dini condition, due to the notion of Dini continuity. Obviously the class $\D$ contains $\phi(s):=\frac{K}{\log^{1+\delta}(c+s^{-1})}$ for constants $K, \delta >0$ and large enough $c\geq \e$ such that $\phi^{2}$ is concave.
\paragraph{Remark 1.2} For any $a\in\A$, we have $\lim_{i\to\infty}a(\lambda_i)=\infty$. The class $\A$ contains a lot of functions. Next, we give a class $\A'$ also containing many functions but the condition in it is more easily to check than the one in $\A$. Letting
\beg{equation*}\beg{split}
\A' = \left\{a\in \B((0,\infty);(0,\infty)), a, \frac{x}{a(x)}\ \  \text{are non-decreasing}, \int_{1}^{\infty}\frac{1}{sa(s)}\d s<\infty \right\},
\end{split}\end{equation*}
we claim $\A'\subset\A$.
\begin{proof} For any $a\in\A'$, $s\in(0,1)$, we have
\beg{equation*}\beg{split}
\sup_{x\geq \frac{1}{s}}\frac{x}{a(x)}\e^{-xs}\leq \sup_{x\geq \frac{1}{s}}\frac{x}{a(\frac{1}{s})}\e^{-xs} \leq \frac{\frac{1}{s}}{a(\frac{1}{s})}\e^{-1}\leq \frac{\frac{1}{s}}{a(\frac{1}{s})}.
\end{split}\end{equation*}
On the other hand,
\beg{equation*}\beg{split}
\sup_{1\wedge\lambda_1\leq x< \frac{1}{s}}\frac{x}{a(x)}\e^{-xs}\leq \sup_{1\wedge\lambda_1\geq x< \frac{1}{s}}\frac{x}{a(x)}\leq \frac{\frac{1}{s}}{a(\frac{1}{s})}.
\end{split}\end{equation*}
So
\beg{equation*}\beg{split}
\int_{0}^{1}\sup_{i\geq 1}\frac{\lambda_i\e^{-\lambda_is}}{a(\lambda_i)}\d s\leq\int_{0}^{1}\sup_{x\in[1\wedge\lambda_1,\infty)}\frac{x}{a(x)}\e^{-xs}\d s\leq \int_{0}^{1}\frac{\frac{1}{s}}{a(\frac{1}{s})}\d s=\int_{1}^{\infty}\frac{1}{sa(s)}\d s<\infty.
\end{split}\end{equation*}
This means $a\in\A$, i.e. $\A'\subset\A$.
\end{proof}
Finally, we give some functions which belong to $\A$.
\beg{enumerate}
\item[{\bf (i)}] $a(x):=x^{\delta}$ for any $\delta\in(0,1]$;

\item[{\bf (ii)}] $a(x):=\log^{1+\delta}(c+x)$ for $\delta>0$ and $c\geq\e^{1+\delta}$;

\item[{\bf (iii)}] $a(x)\geq a_1(x)$, $a_1\in\A'$, for instance, $a(x)=x^{\delta}(\sin x+2)$, $\delta\in(0,1]$.
\end{enumerate}

{\bf (i)} and {\bf (ii)} are in $\A'$. As to {\bf (iii)}, note the fact that if $a\in\A$, then $\tilde{a}\in\B((0,\infty);(0,\infty))$ satisfying $\tilde{a}(x)\geq a(x)$, $x\geq R_0$ for some constant $R_0$ is also in $\A$.

Next, for any $a\in\B((0,\infty);(0,\infty))$, let $\mathbb{H}_{a}=\{x\in\mathbb{H},|a(-A)x|<\infty\}$ equipped the norm $\|x\|_{a}:=|a(-A)x|$, $x\in\mathbb{H}_{a}$. Then $(\mathbb{H}_{a},\|\cdot\|_{a})$ is a Banach space. Note that $\mathbb{H}_{1}=\mathbb{H}$. To obtain the pathwise uniqueness of \eqref{1.1},
we shall need the following condition.
\beg{enumerate}
\item[{\bf (a3)}] For any $(t,x)\in[0,\infty)\times \mathbb{H}$,
\beq\label{1.2}
\lim_{n\to \infty}\|Q(t,x)-Q(t,\pi_{n}x)\|_{\mathrm{HS}}^{2}:=\lim_{n\to \infty}\sum_{k\geq1}|[Q(t,x)-Q(t,\pi_{n}x)]\bar{e}_{k}|^{2}=0.
\end{equation}
Moreover, there exists a function $a\in\A$ such that $b:[0,\infty)\times \mathbb{H}\to \mathbb{H}_{a}$ is measurable and locally bounded, and for any $n\geq1$, there exits $\phi_{n}\in\D$ such that
\beq\label{1.3}
|b(t,x)-b(t,y)|\leq \phi_{n}(|x-y|),\quad t\in[0,n], x,y\in \mathbb{H}, |x|\vee |y|\leq n.
\end{equation}
\end{enumerate}
\paragraph{Remark 1.3} See Remark 3.1 for the reason why we replace \cite[{\bf (a3)}]{W} with the present {\bf (a3)}.

In general, the mild solution (if exists) to \eqref{1.1} can be explosive, so we consider mild solutions with life time.

\beg{defn}  A continuous adapted $\C$-valued process $(X_t)_{t\in[0,\zeta)}$ is called a mild solution to \eqref{1.1} with life time $\zeta$, if $\zeta>0$ is a stopping time such that $\mathbb{P}$-a.s $\limsup_{t\uparrow\zeta}{|X(t)|}=\infty$ holds on $\{\zeta<\infty\}$, and $\mathbb{P}$-a.s
\beg{equation*}\begin{split}
X(t)&= e^{A(t\vee0)} X(t\wedge0)+\int_{0}^{t\vee0} e^{A(t-s)}(b(s,X(s))+B(s,X_{s}))\d s\\
&+\int_{0}^{t\vee0} e^{A(t-s)}Q(s,X(s))\d W(s),\ \  t\in[-r,\zeta)\cap (-\infty,\zeta).
\end{split}\end{equation*}
\end{defn}

The following lemma is a crucial tool in the proof of our results, see \cite[Proposition 7.9]{DZ}.
\beg{lem}\label{L1.1} Let $\{S(t)\}_{t\geq0}$ be a $C_{0}$-contractive semigroup on $\mathbb{H}$. Assume there exists $\alpha\in\left(0,\frac{1}{2}\right)$ and $s>0$ such that
\beq\label{1.4}
\int_{0}^{s}t^{-2\alpha}\|S(t)\|_{\mathrm{HS}}^{2}\d t<\infty.
\end{equation}
Then for every $q\in\left(1,\frac{1}{2\alpha}\right)$, $T>0$, there exists $c_{q}>0$ such that for any
$\L\left(\bar{\mathbb{H}};\mathbb{H}\right)$-valued predictable process $\Phi$, there exists a continuous version of $\int_{0}^{\cdot}S(\cdot-s)\Phi(s)\d W(s)$ such that
\beq\begin{split}\label{1.5}
\mathbb{E}\left[\sup_{t\in[0,T]}\left|\int_{0}^{t}S(t-s)\Phi(s)\d W(s)\right|^{2q}\right]&\leq c_{q}\left[\int_{0}^{T}t^{-2\alpha}\|S(t)\|_{\mathrm{HS}}^{2}\d t\right]^{q}\\
&\times \mathbb{E}\left[\int_{0}^{T}\|\Phi(t)\|^{2q}\d t\right].
\end{split}\end{equation}
\end{lem}
\paragraph{Remark 1.4} Note that {\bf (a1)} implies \eqref{1.4} for $\alpha=\frac{\varepsilon}{2}$ by a simple calculus:
\begin{equation*}\beg{split}
\int_{0}^{s}t^{-2\alpha}\|S(t)\|_{\mathrm{HS}}^{2}\d t&=\sum_{i=1}^{\infty}\int_{0}^{s}t^{-2\alpha}\e^{-2 \lambda_{i}t}\d t\\
&\leq\sum_{i=1}^{\infty}\lambda_i^{2\alpha-1}\int_{0}^{\infty}u^{-2\alpha}\e^{-2 u}\d u\leq C\sum_{i=1}^{\infty}\lambda_i^{2\alpha-1}<\infty.
\end{split}
\end{equation*}
\section{Main results}
\beg{thm}\label{T2.1} Assume {\bf (a1)}, {\bf (a2)} and {\bf (a3)}.
\beg{enumerate}
\item[$(1)$] The equation \eqref{1.1} has a unique mild solution $(X_t)_{t\in[0,\zeta)}$ with life time $\zeta$.
\item[$(2)$] Let $\|Q(t)\|_{\infty}:=\sup_{x\in\mathbb{H}}\|Q(t,x)\|$ be locally bounded in $t\geq 0$. If there exist two positive increasing functions $\Phi, h:[0,\infty)\times[0,\infty)\to(0,\infty)$ such that $\int_{1}^{\infty}\frac{\d s}{\Phi_{t}(s)}=\infty$ for any $t>0$ and
\beq\label{2.1}
\langle B(t,\xi+\eta)+b(t,(\xi+\eta)(0)),\xi(0)\rangle \leq\Phi_{t}\left(\|\xi\|_{\infty}^{2}\right)+h_{t}(\|\eta\|_{\infty}), \ \ \xi, \eta\in\C, t\geq 0,
\end{equation}
then the mild solution is non-explosive.
\end{enumerate}
\end{thm}
For simplicity, we introduced some notations. For any $a\in\B([0,\infty),(0,\infty))$ and $\mathbb{H}_{a}$-valued function $f$ on $[0,T]\times\mathbb{H}$, let
\beg{equation*}
\|f\|_{T,\infty,a}=\sup_{t\in[0,T],x\in\mathbb{H}}|a(-A)f(t,x)|
\end{equation*}
Similarly if $f$ is a operator-valued ( for example $\L(\mathbb{H},\mathbb{H}_{a})$ ) map defined on $[0,T]\times\mathbb{H}$, let
\beg{equation*}
\|f\|_{T,\infty,a}=\sup_{t\in[0,T],x\in\mathbb{H}}\|a(-A)f(t,x)\|.
\end{equation*}
If $a=1$, we omit it.

In order to apply Zvonkin type transformation, we need the following global conditions:
\beg{enumerate}
\item[$\bf{(a2^{'})}$] $B, Q$ satisfy {\bf (a2)}, and there exists a positive increasing function $C_{B,Q}:[0,\infty)\to(0,\infty)$ such that
\beg{equation*}
\|\nabla B\|_{T,\infty}+\sum_{j=0}^{2}\left\|\nabla^{j} Q\right\|_{T,\infty}+\left\|(Q Q^{\ast})^{-1}\right\|_{T,\infty}<C_{B,Q}(T), \ \ T\geq 0.
\end{equation*}
\item[$\bf{(a3^{'})}$] $Q$ satisfies \eqref{1.2}. In addition, for any $T>0$, there exist $a\in\A$ and $\phi\in\D$ such that
    \beq\label{2.2}
\|b\|_{T,\infty,a}<\infty,
\end{equation}
and
\beq\label{2.3}
|b(t,x)-b(t,y)|\leq \phi(|x-y|), \ \ t\in[0,T], x,y\in \mathbb{H}.
\end{equation}
\end{enumerate}

According to Theorem 2.1, under {\bf (a1)}, ${\bf (a2^{'})}$ and ${\bf (a3^{'})}$, the unique mild solution $X^{\xi}_t$ of \eqref{1.1} is non-explosive. Fixing $r<\infty$, the associated Markov semigroup $P_{t}$ of $X_{t}^{\xi}$ is defined as
\beg{equation*}
P_{t}f(\xi)=\mathbb{E}f(X_{t}^{\xi}),\ \ f\in \B_{b}(\C), t\geq 0, \xi\in\C.
\end{equation*}
To derive the Harnack inequalities for $P_{t}$ with $t>r$, we need a stronger condition ${\bf (a3^{''})}$ in stead of ${\bf (a3^{'})}$ as follows:
\begin{enumerate}
\item[$\bf{(a3^{''})}$] $Q$ satisfies \eqref{1.2}. In addition, for any $T>0$,
\begin{equation}\label{2.4}
\left\|(-A)^{\frac{1}{2}}b\right\|_{T,\infty}<\infty,
\end{equation}
and there exists $\phi\in\D$ such that
\begin{equation}\label{2.5}
\left|(-A)^{\frac{1-\varepsilon}{2}}[b(t,x)-b(t,y)]\right|\leq \phi(|x-y|),\quad t\in[0,T], x,y\in \mathbb{H},
\end{equation}
where $\varepsilon$ is in {\bf({a1)}}.
\end{enumerate}
Then we have
\beg{thm}\label{T2.2} Assume {\bf (a1)}, ${\bf (a2^{'})}$, ${\bf (a3^{''})}$. In addition, if $\|B(t)\|_{\infty}:=\sup_{\xi\in\C}|B(t,\xi)|$ is locally bounded in $t\geq 0$, and for any $T>0$, there exists a constant $C(T)>0$ such that
\beq\label{2.6}
\|Q(t,x)-Q(t,y)\|_{\mathrm{HS}}^{2}\leq C(T)|x-y|^{2}, \quad t\in[0,T], x, y\in \mathbb{H}.
\end{equation}
Then for every $T>r$ and positive function $f\in \B_{b}(\C)$,
\beg{enumerate}
\item[(1)] the $\log$-Harnack inequality holds, i.e.
\beq\label{2.7}
P_T\log f(\eta)\leq \log P_T f(\xi)+H(T,\xi,\eta),~\xi,\eta\in\C
\end{equation}
with
$$H(T,\xi,\eta)=C\left(\ff {|\xi(0)-\eta(0)|^2} {T-r}+\|\xi-\eta\|^2_\infty\right)$$
for some constant $C>0$.
\item[(2)]  There exists $K>0$ such that  for any
$p>(1+K)^{2}$,  the Harnack inequality with power
\beq\label{2.8}
P_{T}f(\eta)\leq (P_{T} f^{p}(\xi))^{\frac{1}{p}}\exp{\Psi_{p}(T;\xi,\eta)},\quad \xi,\eta\in\C
\end{equation}
holds, where
\beg{equation*}
\Psi_{p}(T;\xi,\eta)=C(p)\left\{1+\frac{|\xi(0)-\eta(0)|^2}{T-r}+\|\xi-\eta\|_{\infty}^{2}\right\}
\end{equation*}
for a decreasing function $C:\left((1+K)^{2},\infty\right)\to(0,\infty)$.
\end{enumerate}
\end{thm}
The reminder of the paper is organized as follows: In Section 3, we prove the
pathwise uniqueness, in Section 4, combining Section 3 with a
truncating argument, we prove Theorem 2.1, in Section 5, we investigate
the Harnack inequalities for the semigroup  by finite-dimensional approximations. Results in Sections 3 and 5 are derived under some global conditions, and assertions in Section 4 are derived under some local conditions.
\section{Pathwise uniqueness}
In this section, we transform \eqref{1.1} to a regular equation and investigate the pathwise uniqueness of it,   which is equivilant to that of \eqref{1.1}. To this end, firstly, we consider the gradient estimate for the following SPDE \eqref{3.1}, which is crucial in the proof of the regularity of the solution to the equation \eqref{3.3}, see \cite{W} for details. Differently, we need a modified gradient estimate in the present case, and we will give a proof in detail.
\beg{equation}\label{3.1}
\d Z_{s,t}^{x}=A Z_{s,t}^{x}\d t+Q(t,Z_{s,t}^{x})\d W(t), \quad Z_{s,s}^{x}=x, t\geq s\geq 0.
\end{equation}
Then under {\bf (a1)}, ${\bf (a2^{'})}$ with $B=0$, \eqref{3.1} has a unique mild solution $\{Z_{s,t}^{x}\}_{t\geq s}$. Let $P_{s,t}^{0}$ be the associated Markov semigroup.

Firstly, for any $f\in\B_b(\mathbb{H},\mathbb{H})$, $\eta$, $x\in\mathbb{H}$, $0\leq s<t\leq T$, by \cite[(2.12)]{W}, we have
\begin{align*}
\sum_{i=1}^{\infty}\left|\nabla_{\eta} P_{s,t}^{0}\langle f, e_{i}\rangle(x)\right|^{2}\nonumber\leq\sum_{i=1}^{\infty}\frac{c}{t-s}P_{s,t}^{0}|\langle f, e_{i}\rangle|^{2}(x)|\eta|^2=\frac{c}{t-s}P_{s,t}^{0} |f|^{2}(x)|\eta|^2.\nonumber
\end{align*}
Then $\mathbb{H}\ni\nabla_{\eta} P_{s,t}^{0}f(x)\left(:=\sum_{i=1}^{\infty}\left(\nabla_\eta P_{s,t}^0\langle f, e_{i}\rangle(x)\right)e_i\right)$ and
\beg{equation}\label{3.2}
\left|\nabla P_{s,t}^{0}f(x)\right|^{2}\leq\frac{c}{t-s}P_{s,t}^{0} |f|^{2}(x)
\end{equation}
for a constant $c>0$.

Next, by \eqref{3.2}, for any $a\in\B((0,\infty);(0,\infty))$, $f\in \B_{b}(\mathbb{H},\mathbb{H}_a)$ and $x$, $\eta$, $\eta'\in\mathbb{H}$, $0\leq s<t\leq T$, it holds that
\beg{align*}
\sum_{i=1}^\infty a(\lambda_i)^2\left(\nabla_\eta P_{s,t}^0\langle f,e_i\rangle(x)\right)^2&=\sum_{i=1}^\infty \left(\nabla_\eta P_{s,t}^0\langle a(-A)f,e_i\rangle(x)\right)^2\\
&\leq \ff {c} {t-s} P_{s,t}^0|a(-A)f|^2(x)|\eta|^2<\infty.
\end{align*}

Then $\nabla_{\eta}P_{s,t}^0f(x)$ belong to the domain of $a(-A)$ and
\begin{align*}
 a(-A)\nabla_{\eta}P_{s,t}^0f(x)&=\sum_{i=1}^\infty a(\lambda_i)\nabla_{\eta}P_{s,t}^0\langle f,e_i\rangle(x)e_i\\
&=\sum_{i=1}^\infty\nabla_{\eta}P_{s,t}^0\langle a(-A)f,e_i\rangle(x)e_i\\
&=\nabla_{\eta}P_{s,t}^0(a(-A)f)(x).
\end{align*}
Similarly, according to \cite[(2.16)]{W}, it is easy to see that
\beg{equation}\label{3.3}
a(-A)\nabla_{\eta'}\nabla_{\eta}P_{s,t}^0f(x)=\nabla_{\eta'}\nabla_{\eta}P_{s,t}^0(a(-A)f)(x).
\end{equation}
In a word,
\beg{equation}\label{3.4}
a(-A)\nabla^kP_{s,t}^0f=\nabla^kP_{s,t}^0(a(-A)f), ~f\in \B_{b}(\mathbb{H},\mathbb{H}_a), 0\leq s<t\leq T, k=0, 1, 2.
\end{equation}
By \cite{W}, to obtain the pathwise uniqueness of \eqref{1.1}, we need to study the following equation:
\beq\label{3.5}
u(s,x)=\int_{s}^{T} \e^{-\lambda(t-s)}P_{s,t}^{0}(\nabla_{b(t,\cdot)}u(t,\cdot)+b(t,\cdot))(x)\d t, \ \ s\in[0,T].
\end{equation}
The following Lemma is a modified result of \cite[Lemma 2.3]{W}.
\beg{lem}\label{L3.1} Assume {\bf (a1)}, ${\bf (a2^{'})}$ with $B=0$, \eqref{2.2}. Let $T>0$ be fixed. Then there exists a constant $\lambda(T)>0$ such that the following assertions hold.
\beg{enumerate}
\item[(1)] For any $\lambda\geq\lambda(T)$, the equation \eqref{3.5} has a unique solution $u\in C([0,T]; C_{b}^{1}(\mathbb{H}; \mathbb{H}_{a}))$.
\item[(2)] If moreover \eqref{2.3} holds, then we have
\beq\label{3.6}
\lim_{\lambda\to \infty}\|u\|_{T,\infty,a}+\|\nabla u\|_{T,\infty,a}+\left\|\nabla^{2} u\right\|_{T,\infty}=0.
\end{equation}
\end{enumerate}
\end{lem}
\begin{proof}[Proof] (1) Let $\H=C([0,T]; C_{b}^{1}(\mathbb{H}; \mathbb{H}_{a}))$, which is a Banach space under the norm
\begin{equation*}\begin{split}
\|u\|_{\H}:&=\|u\|_{T,\infty,a}+\|\nabla u\|_{T,\infty,a}\\
&=\sup_{t\in[0,T],x\in\mathbb{H}}|a(-A)u(t,x)|+\sup_{t\in[0,T],x\in\mathbb{H}}\|a(-A)\nabla u(t,x)\|,\quad u\in\H.
\end{split}\end{equation*}
For any $u\in\H$, define
\begin{equation*}
(\Gamma u)(s,x)=\int_{s}^{T} \e^{-\lambda(t-s)}P_{s,t}^{0}(\nabla_{b(t,\cdot)}u(t,\cdot)+b(t,\cdot))(x)\d t, \quad s\in[0,T].
\end{equation*}
Then we have $\Gamma\H\subset\H$. In fact, for any $u\in\H$, by ${\bf (a2^{'})}$, \eqref{2.2}, \eqref{3.4} and dominated convergence theorem, it holds that
\begin{equation*}\beg{split}
\|\Gamma u\|_{T,\infty,a}&=\sup_{s\in[0,T],x\in\mathbb{H}}\left|\int_{s}^{T} \e^{-\lambda(t-s)}P_{s,t}^{0}(a(-A)\nabla_{b(t,\cdot)}u(t,\cdot)+a(-A)b(t,\cdot))(x)\d t\right|\\
&\leq \sup_{s\in[0,T]}\int_{s}^{T} \e^{-\lambda(t-s)}(\|b\|_{T,\infty}\|\nabla u\|_{T,\infty,a}+\|b\|_{T,\infty,a})\d t\\
&\leq (\|b\|_{T,\infty}\|\nabla u\|_{T,\infty,a}+\|b\|_{T,\infty,a})\int_{0}^{T} \e^{-\lambda t}\d t\\
&\leq \frac{\|b\|_{T,\infty}\|\nabla u\|_{T,\infty,a}+\|b\|_{T,\infty,a}}{\lambda}<\infty.
\end{split}\end{equation*}
Again by ${\bf (a2^{'})}$, \eqref{2.2}, \eqref{3.4} and dominated convergence theorem, we have
\begin{equation*}\beg{split}
\|\nabla\Gamma u\|_{T,\infty,a}&=\sup_{s\in[0,T],x\in\mathbb{H}, |\eta|\leq 1}\left\|\int_{s}^{T} \e^{-\lambda(t-s)}\nabla_{\eta} P_{s,t}^{0}(a(-A)\nabla_{b(t,\cdot)}u(t,\cdot)+a(-A)b(t,\cdot))(x)\d t\right\|\\
&\leq C\sup_{s\in[0,T]}\int_{s}^{T} \frac{\e^{-\lambda(t-s)}}{\sqrt{t-s}}(\|b\|_{T,\infty}\|\nabla u\|_{T,\infty,a}+\|b\|_{T,\infty,a})\d t\\
&\leq C(\|b\|_{T,\infty}\|\nabla u\|_{T,\infty,a}+\|b\|_{T,\infty,a})\int_{0}^{T} \frac{\e^{-\lambda t}}{\sqrt{t}}\d t\\
&\leq C\frac{\|b\|_{T,\infty}\|\nabla u\|_{T,\infty,a}+\|b\|_{T,\infty,a}}{\sqrt{\lambda}}<\infty.
\end{split}\end{equation*}
So, $\Gamma\H\subset\H$. Next, by the fixed-point theorem, it suffices to show that for large enough $\lambda>0$, $\Gamma$ is contractive on $\H$. To do this, for any $u$, $\tilde{u}\in\H$, it is easy to see that
\begin{equation*}\beg{split}
&\|\Gamma u-\Gamma \tilde{u}\|_{T,\infty,a}\leq \frac{\|b\|_{T,\infty}}{\lambda}\|\nabla u-\nabla \tilde{u}\|_{T,\infty,a},\\
&\|\nabla(\Gamma u-\Gamma \tilde{u})\|_{T,\infty,a}\leq C\frac{\|b\|_{T,\infty}}{\sqrt{\lambda}}\|\nabla u-\nabla \tilde{u}\|_{T,\infty,a}.
\end{split}\end{equation*}
So we can find $\lambda(T)>0$ such that $\Gamma$ is contractive on $\H$ with $\lambda>\lambda(T)$, thus we prove (1).

(2) Combining the proof of (1) and the proof of \cite[Lemma 2.3 (2)]{W}, it is easy to obtain (2). Here to save space, we do not repeat the process.
\end{proof}
The next Lemma gives a regular representation of \eqref{1.1}. See the proof of \cite[Proposition 2.5]{W} for details.
\beg{lem}\label{L3.2} Assume {\bf (a1)}, ${\bf (a2^{'})}$ and ${\bf (a3^{'})}$. Then for any $T>0$, there exists a constant $\lambda(T)>0$ such that for any stopping time $\tau$, any adapted continuous $\C$-valued process $(X_t)_{t\in[0,T\wedge\tau]}$ with $\mathbb{P}$-a.s.
\begin{equation*}\begin{split}
X(t)&= \e^{At} X(0)+\int_{0}^{t} \e^{A(t-s)}(b(s,X(s))+B(s,X_{s}))\d s\\
&+\int_{0}^{t} \e^{A(t-s)}Q(s,X(s))\d W(s), \quad t\in[0,\tau\wedge T],
\end{split}\end{equation*}
and any $\lambda\geq\lambda(T)$, there holds
\beq\label{3.7} \beg{split}
X(t)&= \e^{At} [X(0)+u(0,X(0))]-u(t,X(t))\\
&+\int_{0}^{t}(\lambda-A)\e^{A(t-s)}u(s,X(s))\d s\\
&+\int_{0}^{t}\e^{A(t-s)}[I+\nabla u(s,X(s))]B(s,X_{s})\d s\\
&+\int_{0}^{t} \e^{A(t-s)}[I+\nabla u(s,X(s))]Q(s,X(s))\d W(s), \quad t\in[0,\tau\wedge T],
\end{split}\end{equation}
where $u$ solves \eqref{3.5}, and $\nabla u(s,z)v:=[\nabla_{v} u(s,\cdot)](z)$ for $v$, $z\in\mathbb{H}$.
\end{lem}
\paragraph {Remark 3.1} The second term on the right side of \eqref{3.7} has the same form with the neutral functional SPDE, see \cite{G}. In the case without delay, this can be dealed with by Fubini Theorem, see \cite{W} for details. However, due to the delay, Fubini Theorem is unavailable in the present case. Instead, to prove the pathwise uniqueness, we need a condition like \cite[(H3)]{G}, which can be ensured by Lemma \ref{L3.1} and \eqref{2.2}, see the proof of the following Proposition \ref{P3.3}.

Now, we present a complete proof of the pathwise uniqueness to \eqref{1.1}.
\beg{prp}\label{P3.3} Assume {\bf (a1)}, ${\bf (a2^{'})}$ and ${\bf (a3^{'})}$. Let $\{X_t\}_{t\geq 0},\{Y_t\}_{t\geq 0}$ be two adapted continuous $\C$-valued processes with $X_{0}=Y_{0}=\xi\in\C$. For any $n\geq 1$, let
\beg{equation*}
\tau_{n}^{X}=n\wedge \inf\{t\geq 0:|X(t)|\geq n\}, \ \ \tau_{n}^{Y}=n\wedge \inf\{t\geq 0:|Y(t)|\geq n\}.
\end{equation*}
If $\mathbb{P}$ -a.s. for all $t\in[0,\tau_{n}^{X}\wedge\tau_{n}^{Y}]$, there holds :
\beg{equation*}\beg{split}
&X(t)= \e^{At} \xi(0)+\int_{0}^{t} \e^{A(t-s)}(b(s,X(s))+B(s,X_{s}))\d s+\int_{0}^{t} \e^{A(t-s)}Q(s,X(s))\d W(s),\\
&Y(t)= \e^{At} \xi(0)+\int_{0}^{t} \e^{A(t-s)}(b(s,Y(s))+B(s,Y_{s}))\d s+\int_{0}^{t} \e^{A(t-s)}Q(s,Y(s))\d W(s),
\end{split}\end{equation*}
then $\mathbb{P}$-a.s. $X(t)=Y(t)$, for all $t\in\left[0,\tau_{n}^{X}\wedge\tau_{n}^{Y}\right]$. In particular, $\mathbb{P}$-a.s. $\tau_{n}^{X}=\tau_{n}^{Y}$.
\end{prp}
\beg{proof}[Proof] For any $n\geq 1$, let $\tau_{n}:=\tau_{n}^{X}\wedge\tau_{n}^{Y}$. It suffices to prove that for any $T>0$,
\beq\label{3.8}
\mathbb{E}\sup_{s\in[0,T]}|X(s\wedge\tau_{n})-Y(s\wedge\tau_{n})|^{2p}=0.
\end{equation}
holds for some $p>1$. In the following, we fix $T>0$ and $p>1$. Taking $\lambda$ large enough such that assertions in Lemma \ref{L3.1}, Lemma \ref{L3.2} hold and
\beq\label{3.9}
\frac{5^{4p-1}}{2^{2p+1}}\left(\|a(-A)\nabla u(t,\cdot)\|_{\infty}\int_{0}^{T}\|(-A)[a(-A)]^{-1}\e^{As}\|\d s\right)^{2p}+\|\nabla u(t,\cdot)\|_{\infty}\leq\frac{1}{5}
\end{equation}
for any $t\in[0,T]$.
By \eqref{3.7} for $\tau=\tau_{n}$, we have $\mathbb{P}$-a.s. for any $t\in[0,\tau_{n}\wedge T]$,
\begin{equation*}\beg{split}
&[X(t)+u(t,X(t))]-[Y(t)+u(t,Y(t))]\\
=&\int_{0}^{t}(\lambda-A)\e^{A(t-s)}[u(s,X(s))-u(s,Y(s))]\d s\\
+&\int_{0}^{t}\e^{A(t-s)}\{[I+\nabla u(s,X(s))]B(s,X_{s})-[I+\nabla u(s,Y(s))]B(s,Y_{s})\}\d s\\
+&\int_{0}^{t} \e^{A(t-s)}\{[I+\nabla u(s,X(s))]Q(s,X(s))-[I+\nabla u(s,Y(s))]Q(s,Y(s))\}\d W(s).
\end{split}\end{equation*}

Then \eqref{3.9} yields that
\begin{equation}\label{3.10}\beg{split}
&\mathbb{E}\sup_{t\in[0,q]}|X(t\wedge\tau_{n})-Y(t\wedge\tau_{n})|^{2p}\\
&\leq \frac{5^{4p-1}}{4^{2p}}\mathbb{E}\sup_{t\in[0,q]}\left|\int_{0}^{t\wedge\tau_{n}}(\lambda-A)\e^{A(t-s)}[u(s,X(s))-u(s,Y(s))]\d s\right|^{2p}\\
&+\frac{5^{4p-1}}{4^{2p}}\mathbb{E}\sup_{t\in[0,q]}\left|\int_{0}^{t\wedge\tau_{n}}\e^{A(t-s)}[I+\nabla u(s,X(s))][B(s,X_{s})-B(s,Y_{s})]\d s\right|^{2p}\\
&+\frac{5^{4p-1}}{4^{2p}}\mathbb{E}\sup_{t\in[0,q]}\left|\int_{0}^{t\wedge\tau_{n}}e^{A(t-s)} [\nabla u(s,X(s))-\nabla u(s,Y(s))]B(s,Y_{s})\d s\right|^{2p}\\
&+\frac{5^{4p-1}}{4^{2p}}\mathbb{E}\sup_{t\in[0,q]}\left|\int_{0}^{t\wedge\tau_{n}} \e^{A(t-s)}[\nabla u(s,X(s))-\nabla u(s,Y(s))]Q(s,X(s))\d W(s)\right|^{2p}\\
&+\frac{5^{4p-1}}{4^{2p}}\mathbb{E}\sup_{t\in[0,q]}\left|\int_{0}^{t\wedge\tau_{n}} \e^{A(t-s)}(I+\nabla u(s,Y(s)))[Q(s,X(s))-Q(s,Y(s))]\d W(s)\right|^{2p}\\
&=:I_{1}+I_{2}+I_{3}+I_{4}+I_{5}, \quad q\in[0,T].
\end{split}\end{equation}
Firstly, let
\begin{equation*}
\eta_{q}=\mathbb{E}\sup_{t\in[0,q]}|X(t\wedge\tau_{n})-Y(t\wedge\tau_{n})|^{2p}.
\end{equation*}
Moreover, by \eqref{3.9}, there exists a constant $C(p,\lambda,T)>0$ such that
\beg{equation}\label{3.11}\beg{split}
I_{1}&\leq C(p,\lambda,T)\int_{0}^{q}\eta_{s}\d s\\
&+\left[\frac{5^{4p-1}}{2^{2p+1}}\left(\|a(-A)\nabla u(t,\cdot)\|_{\infty}\int_{0}^{T}\|(-A)[a(-A)]^{-1}\e^{As}\|\d s\right)^{2p}\right]\eta_{q}\\
&\leq C(p,\lambda,T)\int_{0}^{q}\eta_{s}\d s+\frac{1}{5}\eta_{q}.
\end{split}\end{equation}

Since $A$ is negative definite, by \eqref{3.6}, ${\bf (a2^{'})}$, the same initial value of $X$ and $Y$, and H\"{o}lder inequality, it holds that
\beq\label{3.12}\beg{split}
I_{2}&\leq C\mathbb{E}\int_{0}^{q\wedge\tau_{n}}|B(s,X_{s})-B(s,Y_{s})|^{2p}\d s\\
&\leq C_{1}\mathbb{E}\int_{0}^{q}\sup_{t\in[0,s]}|X(t\wedge\tau_{n})-Y(t\wedge\tau_{n})|^{2p}\d s\\
&\leq C_{1}\int_{0}^{q}\eta_{s}\d s
\end{split}\end{equation}
for a constant $C_{1}>0$.
Similarly, Combing \eqref{3.6} and the local boundedness of $B$, we obtain
\begin{equation}\label{3.13}
I_{3}\leq C_{2}\int_{0}^{q}\eta_{s}\d s
\end{equation}
for a constant $C_{2}>0$.

Next, in view of ${\bf (a2^{'})}$, Lemma \ref{L1.1} and Remark 1.4, we have
\beg{equation}\label{3.14}\begin{split}
I_{4}+I_{5}&\leq C_3\mathbb{E}\int_{0}^{q\wedge\tau_{n}}|X(s)-Y(s)|^{2p}\d s=C_3\int_{0}^{q}\eta_{s}\d s,
\end{split}
\end{equation}
for a constant $C_{3}>0$.
Combining \eqref{3.10}, \eqref{3.11}, \eqref{3.12}, \eqref{3.13} and \eqref{3.14}, there exits a constant $C_{0}$ such that
\begin{equation*}
\eta_{l}\leq \frac{1}{5}\eta_{l}+C_{0}\int_{0}^{l}\eta_{q}\d q, \quad l\in[0,T].
\end{equation*}
By Gronwall's inequality, we obtain $\eta_{T}=0$, i.e. \eqref{3.8} holds.
\end{proof}

\section{Proof of Theorem \ref{T2.1}}
\beg{proof}[Proof of Theorem \ref{T2.1}]
(a) We first assume that ${\bf(a1)}$, ${\bf(a2^{'})}$ and ${\bf(a3^{'})}$ hold. Consider the following SPDE on $\mathbb{H}$:
\beg{equation*}
\d Z^{\xi}(t)= A Z^{\xi}(t)\d t +Q(t,Z^{\xi}(t))\d W(t),\ \ Z^{\xi}(0)=\xi(0).
\end{equation*}
It is easy to see that the above equation has a uniqueness non-explosive mild solution:
\beg{equation*}
Z^{\xi}(t)= \e^{At} \xi(0) +\int_{0}^{t} \e^{A(t-s)}Q(s,Z^{\xi}(s))\d W(s),\ \ t\ge 0.
\end{equation*}
Letting $Z^\xi_0=\xi$ (i.e. $Z^\xi(\theta)=\xi(\theta)$ for $\theta \in[-r,0]$), and taking
\beg{equation*}\beg{split}
&W^{\xi}(t)=W(t)-\int_{0}^{t}\psi(s)\d s,\\
&\psi(s)= \big\{Q^{\ast}(QQ^{\ast})^{-1}\big\}(s,Z^{\xi}(s)) \big\{b(s,Z^\xi(s))+ B(s,Z_s^\xi)\big\},\ \ s,t\in[0,T],\end{split}
\end{equation*}
we have
\beg{equation*}\begin{split}
Z^{\xi}(t)&= \e^{At} \xi(0)+\int_{0}^{t} \e^{A(t-s)}B(s,Z_{s}^{\xi})\d s\\
&+\int_{0}^{t} \e^{A(t-s)}b(s,Z^{\xi}(s))\d s+\int_{0}^{t} \e^{A(t-s)}Q(s,Z^{\xi}(s))\d W^{\xi}(s), \ \ t\in[0,T].
\end{split}\end{equation*}
By the local boundedness of $B$, Girsanov theorem implies $\{W^{\xi}(t)\}_{t\in[0,T]}$ is a cylindrical Brownian motion on $\mathbb{\bar{H}}$ under probability $\d\mathbb{Q}^{\xi}=R^{\xi}\d \mathbb{P}$, where
\beg{equation*}
R^{\xi}:=\exp\left[\int_{0}^{T}\big\langle \psi(s),\d W(s)\big\rangle_{\mathbb{\bar{H}}}-\frac{1}{2}\int_{0}^{T}\big|\psi(s)\big|_{\mathbb{\bar{H}}}^{2}\d s\right].
\end{equation*}
Then, under the probability $\Q^\xi$,  $(Z^{\xi}(t),W^{\xi}(t))_{t\in[0,T]}$ is a weak mild solution to \eqref{1.1}.  On the other hand, by Proposition \ref{P3.3},   the pathwise uniqueness holds for the mild solution to \eqref{1.1}. So, by the Yamada-Watanabe principle, the equation \eqref{1.1} has a unique mild solution. Moreover, in this case the solution is non-explosive.

(b) In general, take $\psi\in C_{b}^{\infty}([0,\infty))$ such that $0\leq \psi\leq 1$, $\psi(u)=1$ for $u\in[0,1]$ and $\psi(u)=0$ for $u\in[2,\infty]$. For any $m\geq1$, let
\beg{equation*}\begin{split}
&b^{[m]}(t,z)=b(t\wedge m,z)\psi(|z|/m),\ \ (t,z)\in[0,\infty)\times\mathbb{H},\\
&B^{[m]}(t,\xi)=B(t\wedge m,\xi)\psi(\|\xi\|_{\infty}/m),\ \ (t,\xi)\in[0,\infty)\times\C,\\
&Q^{[m]}(t,z)=Q(t\wedge m,z)\psi(|z|/m),\ \ (t,z)\in[0,\infty)\times\mathbb{H}.
\end{split}\end{equation*}
By {\bf (a2)} and {\bf (a3)}, we know $B^{[m]}$, $Q^{[m]}$ and $b^{[m]}$ satisfy ${\bf (a2^{'})},{\bf (a3^{'})}$. Then by (a), \eqref{1.1} for $B^{[m]}$, $Q^{[m]}$ and $b^{[m]}$ in place of $B$, $Q$, $b$ has a unique mild solution $X^{[m]}(t)$ starting at $X_{0}$ which is non-explosive. Let
\beg{equation*}
\zeta_{0}=0,\ \ \zeta_{m}=m\wedge\inf\{t\geq 0:|X^{[m]}(t)|\geq m\},\ \ m\geq 1.
\end{equation*}
Since $B^{[m]}(s,\xi)=B(s,\xi)$, $Q^{[m]}(s,\xi(0))=Q(s,\xi(0))$ and $b^{[m]}(s,\xi(0))=b(s,\xi(0))$ hold for $s\leq m$, and $\|\xi\|_{\infty}\leq m$,  by Proposition \ref{P3.3}, for any $n$, $m\geq1$, we have $X^{[m]}(t)=X^{[n]}(t)$ for $t\in[0,\zeta_{m}\wedge\zeta_{n}]$. In particular, $\zeta_{m}$ is increasing in $m$. Let $\zeta=\lim_{m\to\infty}\zeta_{m}$ and
\beg{equation*}
X(t)=\sum_{m=1}^{\infty}1_{[\zeta_{m-1},\zeta_{m})}X^{[m]}(t),\ \ t\in[0,\zeta).
\end{equation*}
Then it is easy to see that $X(t)_{t\in[0,\zeta)}$ is a mild solution to \eqref{1.1} with lifetime $\zeta$ and, due to Proposition \ref{P3.3}, the mild solution is unique. So we prove Theorem \ref{T2.1} (1).

(c)Next, we prove the non-explosion.

Let $\|Q\|_{T,\infty}<\infty$ for $T>0$, and let $\Phi, h$ satisfy \eqref{2.1}. Let $X(t)_{t\in[0,\zeta)}$ be the mild solution to \eqref{1.1} with lifetime $\zeta$. Let $M(t)=\int_{0}^{t} e^{A(t-s)}Q(s,X(s))\d W(s)$, $t\in[0,\zeta)$; $M(t)=0$, $t\in[-r,0]$. Then $M_t$ is an adapted continuous process on $\mathbb{H}$ up to the lifetime $\zeta$. It is clear that $Y(t):=X(t)-M(t)$ is the mild solution to the following equation up to $\zeta$,
\beg{equation*}
\d Y(t)= (A Y(t)+b(t,Y(t)+M(t))+B(t,Y_{t}+M_{t}))\d t,\ \ Y_{0}=X_{0}.
\end{equation*}
Then \eqref{2.1} implies that for any $T>0$,
\beg{equation}\label{4.1}\beg{split}
\d |Y(t)|^{2}&\leq 2\langle Y(t), b(t,Y(t)+M(t))+B(t,Y_{t}+M_{t}))\rangle\d t\\
&\leq2\left(\Phi_{\zeta\wedge T}(\|Y_{t}\|_{\infty}^{2})+h_{\zeta\wedge T}(\|M_{t}\|_{\infty})\right)\d t.
\end{split}\end{equation}
Let
\beq\label{4.2}
\Psi_{T}(s)=\int_{1}^{s}\frac{\d r}{2\Phi_{\zeta\wedge T}(r)}, \ \ \alpha_{T}=2\|X_{0}\|_{\infty}^{2}+2\int_{0}^{\zeta\wedge T}h_{\zeta\wedge T}(\|M_{s}\|_{\infty}))\d s.
\end{equation}
It follows from \eqref{4.1} that
\beg{equation}\label{4.3}\begin{split}
\sup_{t\in[0,q]}|Y(t)|^{2}&\leq |X(0)|^{2}+2\int_{0}^{\zeta\wedge T}h_{\zeta\wedge T}(\|M_{s}\|_{\infty}))\d s\\
&+2\int_{0}^{q}\Phi_{\zeta\wedge T}\left(\sup_{t\in[-r,s]}|Y(t)|^{2}\right)\d s, \ \ q\in[0,\zeta\wedge T).
\end{split}\end{equation}
Combining \eqref{4.2} with \eqref{4.3}, we have
\begin{equation}\label{4.4}
\sup_{t\in[-r,q]}|Y(t)|^{2}\leq \alpha_{T}+2\int_{0}^{q}\Phi_{\zeta\wedge T}\left(\sup_{t\in[-r,s]}|Y(t)|^{2}\right)\d s, \ \ q\in[0,\zeta\wedge T).
\end{equation}
Let $Z(s)=\sup_{t\in[-r,s]}|Y(t)|^{2}$, $s\in[0,\zeta\wedge T)$, by Biharis's inequality, \eqref{4.4} implies
\beq\label{4.5}
Z(t)\leq \Psi_{T}^{-1}(\Psi_{T}(\alpha_{T})+t),\ \ t\in[0,\zeta\wedge T).
\end{equation}
Moreover, {\bf (a1)}, $\|Q\|_{T,\infty}<\infty$ and Lemma \ref{L1.1} yield
\beq\label{4.6}
\mathbb{E}\sup_{t\in[0,\zeta\wedge T]}|M(t)|^{2}<\infty.
\end{equation}
So by the definition of $\zeta$ and $Y$, on the set $\{\zeta<\infty\}$, we have $\mathbb{P}$-a.s.
\beq\label{4.7}
\limsup_{t\uparrow\zeta}|Y(t)|=\limsup_{t\uparrow\zeta}|X(t)|=\infty.
\end{equation}
More on the set $\zeta\leq T$, $\mathbb{P}$-a.s. $\alpha_{T}<\infty$. Combining the property of $\Phi$ and \eqref{4.7}, it holds that on the set $\zeta\leq T$, $\mathbb{P}$-a.s.
\beg{equation*}
\infty=\limsup_{t\uparrow\zeta}|Y(t)|^{2}\leq\Psi_{T}^{-1}(\Psi_{T}(\alpha_{T})+T)<\infty.
\end{equation*}
So for any $T>0$, $\mathbb{P}\{\zeta\leq T\}=0$. Note that
\beg{equation*}
 \mathbb{P}\{\zeta< \infty\}=\mathbb{P}\left(\bigcup_{m=1}^{\infty}\{\zeta\leq m\}\right)\leq\sum_{m=1}^{\infty}\mathbb{P}\{\zeta\leq m\}=0,
 \end{equation*}
which implies the solution of \eqref{1.1} is non-explosive.
\end{proof}

\section{Proof of Theorem \ref{T2.2}}

The idea of the proof is to transform \eqref{1.1} into an equation with regular coefficients, so that the Harnack inequalities for the new equation can be derived by coupling by change of measure and finite dimension approximation, see \cite[Theorem 3.4.1, Theorem 4.3.1 and Theorem 4.3.2]{Wbook}. To this end, we use the regularization representation \eqref{3.7}.

In this section, we fix $T>r$. Under {\bf (a1)}, ${\bf (a2^{'})}$, ${\bf (a3^{'})}$ with $a(x)=x^{\frac{1}{2}}$, by Lemma \ref{L3.1} (2) and Lemma \ref{L3.2}, we take large enough $\lambda(T)>0$ such that for any $\lambda\geq\lambda(T)$, Lemma \ref{L3.2} holds and the unique solution $u$ to \eqref{3.5} satisfies:
\beq\label{5.1}\begin{split}
\|\nabla^{2} u\|_{T,\infty}\leq\frac{1}{8}, \ \ \|(-A)^{\frac{1}{2}}\nabla u\|_{T,\infty}\leq\frac{\sqrt{\lambda_1}}{8}.
\end{split}\end{equation}
To treat the delay part, define $u(s,\cdot)=u(0,\cdot)$ for $s\in[-r,0]$. Let $\theta(t,x)=x+u(t,x)$, $(t,x)\in[-r,T]\times\mathbb{H}$.
By \eqref{5.1}, $\{\theta(t,\cdot)\}_{t\in[-r,T]}$ is a family of diffeomorphisms on $\mathbb{H}$. For simplicity, we write $\theta^{-1}(t,x)=[\theta^{-1}(t,\cdot)](x)$, $(t,x)\in[-r,T]\times\mathbb{H}$.
By \eqref{5.1}, we have
\beq\label{5.2}
\frac{7}{8}\leq\|\nabla \theta(t,x)\|\leq\frac{9}{8},\quad \frac{8}{9}\leq\|\nabla \theta^{-1}(t,x)\|\leq\frac{8}{7},\ \ (t,x)\in[-r,T]\times\mathbb{H}.
\end{equation}
where $\nabla \theta(t,x):=[\nabla \theta(t,\cdot)](x)$ and $\nabla \theta^{-1}(t,x):=[\nabla \theta^{-1}(t,\cdot)](x)$.

On the other hand, for any $t\in[0,T]$, define $\theta_t:\C\to\C$ as
\beg{equation}\label{5.3}\begin{split}
(\theta_{t}(\xi))(s)=\theta(t+s,\xi(s)),\ \ \xi\in\C, s\in[-r,0].
\end{split}\end{equation}
Then $\{\theta_{t}\}_{t\in[0,T]}$ is a family of diffeomorphisms on $\C$. Moreover, it is easy to see that for any $t\in[0,T]$,
\beg{equation}\label{5.4}\begin{split}
(\theta_{t}^{-1}(\xi))(s)=\theta^{-1}(t+s,\xi(s)),\ \ \xi\in\C, s\in[-r,0].
\end{split}\end{equation}

Furthermore, it follows from \eqref{5.3} and \eqref{5.2} that
\beg{equation}\begin{split}\label{5.5}
\|(\nabla \theta_{t})(\xi)\|:&=\limsup_{\eta\to \xi}\frac{\sup_{s\in[-r,0]}|\theta(t+s,\eta(s))-\theta(t+s,\xi(s))|}{\|\eta-\xi\|_{\infty}}\\
&\leq\frac{9}{8}, \quad t\in[0,T],\xi\in\C.
\end{split}\end{equation}
Similarly, we have
\beg{equation}\begin{split}\label{5.6}
\|(\nabla \theta_{t}^{-1})(\xi)\|\leq\frac{8}{7}, \quad t\in[0,T],\xi\in\C.
\end{split}\end{equation}

Now, letting $\{X^{\xi}(t)\}_{t\in[-r,T]}$ solve \eqref{1.1} with $X_{0}^{\xi}=\xi\in\C$, by \eqref{3.7}, for any $\lambda\geq \lambda(T)$, $\{Y^{\xi}(t)=\theta(t,X^{\xi}(t))\}_{t\in[-r,T]}$ with $Y^{\xi}_{t}=\theta_{t}(X^{\xi}_{t})$ satisfies
\beq\label{5.7}\beg{split}
Y^{\xi}(t)&= e^{At} Y^{\xi}(0)+\int_{0}^{t}e^{A(t-s)}(\lambda-A)u\left(s,\theta^{-1}\left(s, Y^{\xi}(s)\right)\right)\d s\\
&+\int_{0}^{t}e^{A(t-s)}\nabla \theta\left(s, \theta^{-1}\left(s,Y^{\xi}(s)\right)\right)B\left(s,\theta_{s}^{-1}\left(Y^{\xi}_{s}\right)\right)\d s\\
&+\int_{0}^{t} e^{A(t-s)}\nabla \theta\left(s,\theta^{-1}\left(s,Y^{\xi}(s)\right)\right)Q\left(s,\theta^{-1}\left(s,Y^{\xi}(s)\right)\right)\d W(s),\quad t\in[0,T].
\end{split}\end{equation}
Let
\beg{equation}\label{5.8}
\bar{b}(t, x)=(\lambda-A)u\left(t,\theta^{-1}(t,x)\right),\quad t\in[0,T], x\in\mathbb{H}.
\end{equation}
\beg{equation}\label{5.9}
\bar{B}(t,\xi)=\nabla\theta\left(t,\theta^{-1}(t,\xi(0))\right)B\left(t,\theta_{t}^{-1}(\xi)\right),\quad t\in[0,T], \xi\in\C.
\end{equation}
\beg{equation}\label{5.10}
\bar{Q}(t,x)=\nabla\theta\left(t,\theta^{-1}(t,x)\right)Q\left(t,\theta^{-1}(t,x)\right), \quad t\in[0,T], x\in\mathbb{H}.
\end{equation}
Then for any $\lambda\geq\lambda(T)$, $\{\bar{X}^{\xi}(t):=Y^{\theta_{0}^{-1}(\xi)}(t)\}_{t\in[-r,T]}$ is a mild solution to the equation
\beq\label{5.11}
\d \bar{X}^{\xi}(t)= \left[A \bar{X}^{\xi}(t)+\bar{b}\left(t,\bar{X}^{\xi}(t)\right)+\bar{B}\left(t,\bar{X}^{\xi}_{t}\right)\right]\d t+\bar{Q}\left(t,\bar{X}^{\xi}(t)\right)\d W(t),\quad \bar{X}^{\xi}_{0}=\xi.
\end{equation}
Define
\beg{equation*}
\bar{P}_{t}f(\xi)=\mathbb{E}f(\bar{X}_{t}^{\xi}),\quad t\in[0,T], f\in \B_{b}(\C).
\end{equation*}
Then it is easy to see that
\beq\label{5.12}\beg{split}
P_{t}f(\xi):&=\mathbb{E}f(X_{t}^{\xi})=\mathbb{E}(f\circ\theta_{t}^{-1})(Y^{\xi}_{t})=
\mathbb{E}(f\circ\theta_{t}^{-1})(\bar{X}^{\theta_{0}(\xi)}_{t})\\
&=\bar{P}_{t}(f\circ\theta_{t}^{-1})(\theta_{0}(\xi)), \quad\xi\in\C,t\in[0,T], f\in \B_{b}(\C).
\end{split}\end{equation}
We first study the Harnack inequalities for $\bar{P}_{t}$.

To apply the method of coupling by change of measure, we will use the finite dimension approximation argument. More precisely, let $\{\bar{X}^{n,\xi}(t)\}_{t\in[-r,T]}$ solves the finite-dimensional equation on $\mathbb{H}_{n}:=\mathrm{span}\{e_1,\cdots,e_n\}$ ($n\geq 1$):
\beq\label{5.13}\beg{split}
\d \bar{X}^{(n,\xi)}(t)&= \left[A \bar{X}^{(n,\xi)}(t)+\bar{b}^{n}\left(t,\bar{X}^{(n,\xi)}(t)\right)+\bar{B}^{n} \left(t,\bar{X}^{(n,\xi)}_{t}\right)\right]\d t\\
&+\bar{Q}^{n}\left(t,\bar{X}^{(n,\xi)}(t)\right)\d W(t), \quad\bar{X}^{(n,\xi)}_{0}=\xi\in\C(\mathbb{H}_{n}),
\end{split}\end{equation}
where $\bar{b}^{n}=\pi_{n}\bar{b}$, $\bar{B}^{n}=\pi_{n}\bar{B}$, $\bar{Q}^{n}=\pi_{n}\bar{Q}$. Firstly, we prove that there exists a constant $\tilde{\lambda}(T)\geq\lambda(T)$ such that for any $\lambda\geq \tilde{\lambda}(T)$,
\begin{equation}\label{5.14}
\lim_{n\to \infty}\mathbb{E}\left\|\bar{X}^{\xi}_{t}-\bar{X}^{(n,\pi_{n}\xi)}_{t}\right\|_{\infty}^{2}=0, \ \ t\in[0,T].
\end{equation}
\beg{lem}\label{L5.1} Assume {\bf (a1)}, ${\bf (a2^{'})}$ and ${\bf (a3^{''})}$ with \eqref{2.3} in place of \eqref{2.5}. If in addition $\|B\|_{T,\infty}<\infty$,
then there exists a constant $\tilde{\lambda}(T)\geq\lambda(T)$ such that for any $\lambda\geq \tilde{\lambda}(T)$, \eqref{5.14} holds.
\end{lem}
\beg{proof}[Proof] For simplicity, we omit $\xi$ and $\pi_{n}\xi$ from the subscripts, i.e. we write $(\bar{X}_{t}, \bar{X}^{(n)}_{t})$ instead of $(\bar{X}^{\xi}_{t}, \bar{X}^{(n,\pi_{n}\xi)}_{t})$.
By Jensen inequality, it suffices to prove there exists a constant $\tilde{\lambda}(T)\geq\lambda(T)$ such that for any $\lambda\geq \tilde{\lambda}(T)$,
\begin{equation}\label{5.15}
\lim_{n\to \infty}\mathbb{E}\left\|\bar{X}^{\xi}_{t}-\bar{X}^{(n,\pi_{n}\xi)}_{t}\right\|_{\infty}^{4}=0, \ \ t\in[0,T].
\end{equation}
For any $t\in[0,T]$, let
\beg{equation*}
\beta_{n}(t)=\mathbb{E}\sup_{s\in[-r,t]}|\bar{X}(s)-\bar{X}^{(n)}(s)|^{4},
\ \ \beta(t)=\limsup_{n\to\infty}\beta_{n}(t).
\end{equation*}
Obviously, {\bf (a1)}, ${\bf (a2^{'})}$ and ${\bf (a3^{''})}$ imply
\beg{equation*}
\mathbb{E}\sup_{t\in[0,T], n\geq 1}(\|\bar{X}_{t}\|_{\infty}^{4}+\|\bar{X}^{(n)}_{t}\|_{\infty}^{4})<\infty,
\end{equation*}
so that $\beta(t)<\infty$.

Combining \eqref{5.11} with \eqref{5.13}, it holds that
\beq\label{5.16}\beg{split}
\beta_{n}(t)&\leq 686\|\xi-\pi_{n}\xi\|_{\infty}^{4}\\
&+343\mathbb{E}\sup_{q\in[0,t]}\left|\int_{0}^{q} e^{A(q-s)}[\bar{b}^{n}(s,\bar{X}(s))-\bar{b}^{n}(s,\bar{X}^{(n)}(s))]\d s\right|^{4}\\
&+343\mathbb{E}\sup_{q\in[0,t]}\left|\int_{0}^{q} e^{A(q-s)}[\bar{b}(s,\bar{X}(s))-\bar{b}^{n}(s,\bar{X}(s))]\d s\right|^{4}\\
&+343\mathbb{E}\sup_{q\in[0,t]}\left|\int_{0}^{q} e^{A(q-s)}[\bar{B}^{n}(s,\bar{X}_{s})
-\bar{B}^{n}(s,\bar{X}^{(n)}_{s})]\d s\right|^{4}\\
&+343\mathbb{E}\sup_{q\in[0,t]}\left|\int_{0}^{q} e^{A(q-s)}[\bar{B}(s,\bar{X}_{s})
-\bar{B}^{n}(s,\bar{X}_{s})]\d s\right|^{4}\\
&+343\mathbb{E}\sup_{q\in[0,t]}\left|\int_{0}^{q} e^{A(q-s)}[\bar{Q}(s,\bar{X}(s))-\bar{Q}^{n}(s,\bar{X}(s))]\d W(s)\right|^{4}\\
&+343\mathbb{E}\sup_{q\in[0,t]}\left|\int_{0}^{q} e^{A(q-s)}[\bar{Q}^{n}(s,\bar{X}(s))-\bar{Q}^{n}(s,\bar{X}^{(n)}(s))]\d W(s)\right|^{4}\\
&=:\Gamma_{1}+\Gamma_{2}+\Gamma_{3}+\Gamma_{4}+\Gamma_{5}+\Gamma_{6}+\Gamma_{7}.
\end{split}\end{equation}
Firstly, for any $\xi\in\C$, $n\geq 1$, by the definition of $\pi_{n}$, we have $\|\xi-\pi_{n}\xi\|_{\infty}<\|\xi\|_\infty$ and $|(\xi-\pi_{n}\xi)(s_1)-(\xi-\pi_{n}\xi)(s_2)|\leq |\xi(s_1)-\xi(s_2)|$ for any $s_1$, $s_2\in[-r,0]$. Since for any $s\in[-r,0]$, $|(\xi-\pi_{n}\xi)(s)|\to\infty$ as $n\to\infty$, it follows from Arzela-Ascoli Theorem that
\beg{equation}\label{5.17}
\lim_{n\to\infty}\Gamma_{1}=0.
\end{equation}
Similarly to the estimate $I_{1}$ in Proposition \ref{P3.3}, there exists a constant $\tilde{\lambda}(T)\geq\lambda(T)$ such that for any $\lambda\geq \tilde{\lambda}(T)$,
\beg{equation}\label{5.18}\begin{split}
\Gamma_{2}&\leq C(\lambda,T)\int_{0}^{t}\beta_n(s)\d s+\frac{1}{5}\mathbb{E}\beta_n(t).
\end{split}\end{equation}

Next, by H\"{o}lder inequality and ${\bf (a3^{''})}$, for any $\delta\in(0,2)$, it holds that
\beg{equation*}\begin{split}
\Gamma_{3}&\leq C\e^{\lambda T} \left\{\int_{0}^{t}\left\|(\lambda-A)^{\frac{1}{2}}e^{-(\lambda-A)s}\right\|^{\delta}\d s\right\}^{\frac{4}{\delta}}\\ &\times\mathbb{E}\left\{\int_{0}^{t}\left|(\lambda-A)^\frac{1}{2}(u-\pi_{n}u)\left(s,\theta^{-1}\left(s,\bar{X}(s)\right)\right) \right|^{\frac{\delta}{\delta-1}}\d s\right\}^{\frac{4(\delta-1)}{\delta}}\\
&\leq C(\lambda,T,\delta)\mathbb{E}\left\{\int_{0}^{t}\left|(\lambda-A)^\frac{1}{2}(u-\pi_{n}u)\left(s,\theta^{-1}\left(s,\bar{X}(s)\right)\right) \right|^{\frac{\delta}{\delta-1}}\d s\right\}^{\frac{4(\delta-1)}{\delta}}.
\end{split}\end{equation*}
Combing the definition of $\pi_{n}$ and Lemma \ref{3.1} (2) for $a(x)=x^{\frac{1}{2}}$, it follows from dominated convergence theorem that
\beg{equation}\label{5.19}
\lim_{n\to\infty}\Gamma_{3}=0.
\end{equation}
Moreover, ${\bf (a2^{'})}$, H\"{o}lder inequality and the boundedness of $B$ yield that
\beg{equation}\begin{split}\label{5.20}
\Gamma_{4}&\leq C(\lambda, T)\int_{0}^{t}\beta_{n}(s)\d s.
\end{split}\end{equation}
Again using the boundedness of $B$, H\"{o}lder inequality and dominated convergence theorem, we obtain
\beg{equation}\label{5.21}
\lim_{n\to\infty}\Gamma_{5}=0.
\end{equation}
Furthermore, combining ${\bf (a2^{'})}$ with Lemma \ref{L1.1}, applying dominated convergence theorem, it is easy to see that
\beg{equation}\label{5.22}
\lim_{n\to\infty}\Gamma_{6}=0.
\end{equation}
Finally, combining ${\bf (a2^{'})}$ with Lemma \ref{L1.1}, we have
\beg{equation}\begin{split}\label{5.23}
\Gamma_{7}&\leq C(\lambda, T)\int_{0}^{t}\beta_{n}(s)\d s.
\end{split}\end{equation}
Combining \eqref{5.16}-\eqref{5.23}, applying dominated convergence theorem, it holds that
\beg{equation*}
\beta(t)\leq C\int_{0}^{t}\beta(s)\d s.
\end{equation*}
Since $\beta(t)<\infty$, Gronwall inequality yields $\beta(t)=0$, $t\in[0,T]$, which implies \eqref{5.14}.
\end{proof}
\beg{lem}\label{L5.2} Assume {\bf (a1)}, ${\bf (a2^{'})}$ with $B=0$, \eqref{2.4} and \eqref{2.5}.
Then for any $\lambda\geq \lambda(T)$, there exists a constant $C(T)>0$ such that
\beg{equation}\label{5.24}
\|\nabla u(t,x)-\nabla u(t,y)\|_{\mathrm{HS}}\leq C(T)|x-y|, \quad x, y\in\mathbb{H}, t\in[0,T].
\end{equation}
\end{lem}
\begin{proof}
In order to prove \eqref{5.24}, by {\bf (a1)}, it suffices to prove
\beg{equation}\label{5.25}
\left\|(-A)^{\frac{1-\varepsilon}{2}}[\nabla u(t,x)-\nabla u(t,y)]\right\|\leq C(T)|x-y|, \quad x, y\in\mathbb{H}, t\in[0,T].
\end{equation}
In fact, if \eqref{5.25} holds, then
\beg{equation*}\begin{split}
\|\nabla u(t,x)-\nabla u(t,y)\|_{\mathrm{HS}}^{2}
&=\left\|(-A)^{\frac{\varepsilon-1}{2}}(-A)^{\frac{1-\varepsilon}{2}}[\nabla u(t,x)-\nabla u(t,y)]\right\|_{\mathrm{HS}}^{2}\\
&\leq C(T)\left\|(-A)^{\frac{\varepsilon-1}{2}}\right\|_{\mathrm{HS}}^{2}|x-y|^{2}, \quad x, y\in\mathbb{H}, t\in[0,T].
\end{split}\end{equation*}
Define
\beg{equation*}\begin{split}
(R^{\lambda}_{s,t}f)(x)=\int_{s}^{t}\e^{-(q-s)\lambda}[P_{s,q}^{0}f(q,\cdot)](x), \ \ x\in\mathbb{H}, \lambda\geq 0, t\geq s\geq 0, f\in\B_{b}([0,\infty)\times \mathbb{H};\mathbb{H}).
\end{split}\end{equation*}
Firstly, by \eqref{2.4} and \eqref{5.1}, we have
\begin{equation}\label{5.26}
\|(-A)^{\frac{1}{2}}(\nabla_{b}u+b)\|_{T,\infty}<\infty
\end{equation}
for any $\lambda\geq\lambda(T)$.

Then it follows from \eqref{3.4}, \eqref{3.5}, \eqref{5.26}, \cite[Lemma 2.2 (1)]{W} and dominated convergence theorem that for any $\lambda\geq\lambda(T)$,
\beg{equation}\begin{split}\label{5.27}
&\left\|(-A)^{\frac{1}{2}}[\nabla u(t,x)-\nabla u(t,y)]\right\|\\
=&\left\|\nabla\left(R^{\lambda}_{t,T}\left((-A)^{\frac{1}{2}}(\nabla_{b}u+b)\right)\right)(x)- \nabla\left(R^{\lambda}_{t,T}\left((-A)^{\frac{1}{2}}(\nabla_{b}u+b)\right)\right)(y)\right\|\\
\leq& C|x-y|\log\left(\e+\frac{1}{|x-y|}\right) \quad x,y\in\mathbb{H}, t\in[0,T]
\end{split}\end{equation}
holds for some constant $C>0$. Combining this with \eqref{5.26} and \eqref{2.5}, it is easy to see that for any $\lambda\geq\lambda(T)$,
\beg{equation}\label{5.28}
\left|(-A)^{\frac{1-\varepsilon}{2}}(\nabla_{b}u+b)(t,x)-(-A)^{\frac{1-\varepsilon}{2}}(\nabla_{b}u+b)(t,y)\right|\leq \tilde{\phi}(|x-y|), \ \ t\in[0,T], x,y\in\mathbb{H},
\end{equation}
where $\tilde{\phi}(s)=c\sqrt{\phi^2(s)+s}$ with a constant $c>0$.

Finally, by \eqref{3.4}, \eqref{3.5}, \eqref{5.28}, dominated convergence theorem and \cite[Lemma 2.2 (3)]{W}, for any $\lambda\geq \lambda(T)$, we conclude that
\beg{equation}\begin{split}\label{5.29}
&\left\|(-A)^{\frac{1-\varepsilon}{2}}[\nabla u(t,x)-\nabla u(t,y)]\right\|\\
=&\left\|\nabla\left(R^{\lambda}_{t,T}\left((-A)^{\frac{1-\varepsilon}{2}}(\nabla_{b}u+b)\right)\right)(x)- \nabla\left(R^{\lambda}_{t,T}\left((-A)^{\frac{1-\varepsilon}{2}}(\nabla_{b}u+b)\right)\right)(y)\right\|\\
\leq& C(T)|x-y|, \quad x,y\in\mathbb{H}, t\in[0,T].
\end{split}\end{equation}
for a constant $C(T)>0$. Thus \eqref{5.25} holds, and we complete the proof.
\end{proof}

\beg{lem}\label{L5.3} Assume {\bf (a1)}, ${\bf (a2^{'})}$, \eqref{2.4} and \eqref{2.5}. If in addition $\|B\|_{T,\infty}<\infty$ and moreover
\beq\label{5.30}
\|Q(t,x)-Q(t,y)\|_{\mathrm{HS}}^{2}\leq C(T)|x-y|^{2},\ \ t\in[0,T],x, y \in \mathbb{H},
\end{equation}
where $C(T)$ is a positive constant. Then for any $\lambda\geq\lambda(T)$, there exists $K_{1}\geq 0$, $K_{2}\geq 0$, $K_{3}>0$ and $K_{4}\in \mathbb{R}$ ($K_{1}$, $K_{2}$, $K_{3}$, $K_{4}$ only depend on T) such that
\beg{equation}\label{5.31}
\left|(\bar{Q}^{\ast}(\bar{Q}\bar{Q}^{\ast})^{-1})(t,\eta(0))\{\bar{B}(t,\xi)-\bar{B}(t,\eta)\}\right|_{\bar{\mathbb{H}}}\leq K_{1}\|\xi-\eta\|_{\infty};
\end{equation}
\beg{equation}\label{5.32}
\left\|\bar{Q}(t,x)-\bar{Q}(t,y)\right\|\leq K_{2}(1\wedge|(x-y|);
\end{equation}
\beg{equation}\label{5.33}
\left\|(\bar{Q}^{\ast}(\bar{Q}\bar{Q}^{\ast})^{-1})(t,x)\right\|\leq K_{3};
\end{equation}
\beg{equation}\label{5.34}
\left\|\bar{Q}(t,x)-\bar{Q}(t,y)\right\|_{HS}^{2}+2\left\langle x-y, Ax-Ay+\bar{b}(t,x)-\bar{b}(t,y)\right\rangle\leq K_{4}|x-y|^{2};
\end{equation}

hold for $t\in[0,T]$, $\xi,\eta\in\C$, and $x,y\in\mathbb{H}$.
\end{lem}
\beg{proof}[Proof] Fix $\lambda\geq \lambda(T)$.

(a) Since $\nabla\theta(t,\cdot)=I+\nabla u(t,\cdot)$, $t\in[0,T]$, \eqref{5.1} yields that for any $\lambda\geq\lambda(T)$ and $(t,x)\in[0,T]\times\mathbb{H}$, $\nabla\theta (t,x)$, $(\nabla\theta(t,x))^{\ast}\in\L(\mathbb{H},\mathbb{H})$ are invertible. Then from \eqref{5.10},
\beg{equation}\beg{split}\label{5.35}
\left(\bar{Q}^{\ast}(\bar{Q}\bar{Q}^{\ast})^{-1}\right)(t,x)= \left(Q^{\ast}(QQ^{\ast})^{-1}\right) \left(t,\theta^{-1}(t,x)\right)\left[\nabla\theta\left(t,\theta^{-1}(t,x)\right)\right]^{-1}
\end{split}
\end{equation}
From ${\bf (a2^{'})}$, \eqref{5.33} holds with $K_{3}=\frac{8}{7}C_{B,Q}^{2}(T)$.

(b) Due to (a), in order to prove \eqref{5.31}, we only need to estimate $\left|\bar{B}(t,\xi)-\bar{B}(t,\eta)\right|_{\mathbb{H}}$. From \eqref{5.9}, ${\bf{(a2^{'})}}$, \eqref{5.1}, \eqref{5.2}, we have
\beg{equation}\beg{split}\label{5.36}
&|\bar{B}(t,\xi)-\bar{B}(t,\eta)|\\
=&\left|\nabla\theta(t,\theta^{-1}(t,\xi(0)))B(t,\theta_{t}^{-1}(\xi)) -\nabla\theta(t,\theta^{-1}(t,\eta(0)))B(t,\theta_{t}^{-1}(\eta))\right|\\
\leq&\left|\nabla\theta(t,\theta^{-1}(t,\xi(0)))B(t,\theta_{t}^{-1}(\xi)) -\nabla\theta(t,\theta^{-1}(t,\xi(0)))B(t,\theta_{t}^{-1}(\eta))\right|\\
+&\left|\nabla\theta(t,\theta^{-1}(t,\xi(0)))B(t,\theta_{t}^{-1}(\eta)) -\nabla\theta(t,\theta^{-1}(t,\eta(0)))B(t,\theta_{t}^{-1}(\eta))\right|\\
\leq& \sup_{(t,x)\in[0,T]\times\mathbb{H}}\|\nabla\theta(t,x)\|\|\nabla B\|_{T,\infty}\sup_{t\in[0,T],\xi\in\C}\left\|\nabla\theta_{t}^{-1}(\xi)\right\|\|\xi-\eta\|_{\infty}\\
+&\sup_{(t,x)\in[0,T]\times\mathbb{H}}\|\nabla^{2}\theta(t,x) \|\sup_{(t,x)\in[0,T]\times\mathbb{H}}\|\nabla\theta^{-1}(t,x)\|\|B\|_{T,\infty}|\xi(0)-\eta(0)|\\
\leq& K\|\xi-\eta\|_{\infty}, \ \ K>0.
\end{split}\end{equation}
Combining \eqref{5.36} with \eqref{5.33}, we prove \eqref{5.31}.

(c) Similarly, from\eqref{5.10}, again using  ${\bf{(a2^{'})}}$, \eqref{5.1}, \eqref{5.2}, we arrive at
\beg{equation}\beg{split}\label{5.37}
&\left\|\bar{Q}(t,x)-\bar{Q}(t,y)\right\|\\
=&\left\|\nabla\theta(t,\theta^{-1}(t,x)) Q(t,\theta^{-1}(t,x))-\nabla\theta(t,\theta^{-1}(t,y)) Q(t,\theta^{-1}(t,y))\right\|\\
\leq&\left\|\nabla\theta(t,\theta^{-1}(t,x)) Q(t,\theta^{-1}(t,x))-\nabla\theta(t,\theta^{-1}(t,x)) Q(t,\theta^{-1}(t,y))\right\|\\
+&\left\|\nabla\theta(t,\theta^{-1}(t,x)) Q(t,\theta^{-1}(t,y))-\nabla\theta(t,\theta^{-1}(t,y)) Q(t,\theta^{-1}(t,y))\right\|\\
\leq& \sup_{(t,x)\in[0,T]\times\mathbb{H}}\|\nabla\theta(t,x)\|\|\nabla Q\|_{T,\infty}\sup_{(t,x)\in[0,T]\times\mathbb{H}}\|\nabla\theta^{-1}(t,x)\||x-y|\\
+&\sup_{(t,x)\in[0,T]\times\mathbb{H}}\|\nabla^{2}\theta(t,x)\| \|Q\|_{T,\infty}\sup_{(t,x)\in[0,T]\times\mathbb{H}}\|\nabla\theta^{-1}(t,x)\||x-y|\\
\leq& K^{'}|x-y|, \ \ K^{'}>0,
\end{split}
\end{equation}
and
\beg{equation}\beg{split}\label{5.38}
&\|\bar{Q}(t,x)-\bar{Q}(t,y)\|\leq 2\sup_{(t,x)\in[0,T]\times\mathbb{H}}\|\nabla\theta(t,x)\|\| Q\|_{T,\infty}\leq K^{''},\quad K^{''}>0.
\end{split}\end{equation}
Then \eqref{5.37} and \eqref{5.38} yield \eqref{5.32}.

(d) Finally, from\eqref{5.10}, applying  ${\bf{(a2^{'})}}$, \eqref{5.1}, \eqref{5.2}, \eqref{5.30} and Lemma \ref{L5.2}, we obtain
\beg{equation}\beg{split}\label{5.39}
&\left\|\bar{Q}(t,x)-\bar{Q}(t,y)\right\|_{\mathrm{HS}}\\
=&\left\|\nabla\theta(t,\theta^{-1}(t,x)) Q(t,\theta^{-1}(t,x))-\nabla\theta(t,\theta^{-1}(t,y)) Q(t,\theta^{-1}(t,y))\right\|_{\mathrm{HS}}\\
\leq&\left\|\nabla\theta(t,\theta^{-1}(t,x)) Q(t,\theta^{-1}(t,x))-\nabla\theta(t,\theta^{-1}(t,x)) Q(t,\theta^{-1}(t,y))\right\|_{\mathrm{HS}}\\
+&\left\|\nabla\theta(t,\theta^{-1}(t,x)) Q(t,\theta^{-1}(t,y))-\nabla\theta(t,\theta^{-1}(t,y)) Q(t,\theta^{-1}(t,y))\right\|_{\mathrm{HS}}\\
\leq& \sup_{(t,x)\in[0,T]\times\mathbb{H}}\|\nabla\theta(t,x)\| \left\|Q(t,\theta^{-1}(t,x))-Q(t,\theta^{-1}(t,x)\right\|_{\mathrm{HS}} \\
+ &\|Q\|_{T,\infty}\left\|\nabla u(t,\theta^{-1}(t,x))-\nabla u(t,\theta^{-1}(t,x))\right\|_{\mathrm{HS}}\\
\leq& C(T)\left[\sup_{t\in[0,T]\times\mathbb{H}}\|\nabla\theta(t,x)\|+ \|Q\|_{T,\infty}\right]\sup_{t\in[0,T]\times\mathbb{H}}\left\|\nabla\theta^{-1}(t,x)\right\||x-y|\\
\leq& K_{0}|x-y|, \ \ K_{0}>0.
\end{split}\end{equation}
Moreover, \eqref{5.8} \eqref{5.1}, \eqref{5.2}, and ${\bf{(a3^{''})}}$ yield that
\beg{equation}\beg{split}\label{5.40}
&\langle A(x-y), x-y\rangle+\left\langle(-A)[u(t,\theta^{-1}(t,x))-u(t,\theta^{-1}(t,y))], x-y\right\rangle\\
=&-\left|(-A)^{\frac{1}{2}}(x-y)\right|^{2}+\left\langle (-A)^{\frac{1}{2}}\left[u\left(t,\theta^{-1}(t,x)\right)-u\left(t,\theta^{-1}(t,y)\right)\right],(-A)^{\frac{1}{2}} (x-y)\right\rangle\\
\leq&-|(-A)^{\frac{1}{2}}(x-y)|^{2}+c|x-y|^{2}+\frac{1}{2}\left|(-A)^{\frac{1}{2}}(x-y)\right|^{2}\\
\leq& c|x-y|^{2}
\end{split}\end{equation}
for a constant $c>0$. Combining \eqref{5.39} with \eqref{5.40}, we obtain \eqref{5.34}.
\end{proof}
\begin{proof}[Proof of Theorem \ref{T2.2}] Combining Lemma \ref{L5.1} with Lemma \ref{L5.3}, for any $\lambda\geq\tilde{\lambda}(T)$, using \cite[Theorem 3.4.1, Theorem 4.3.1 and Theorem 4.3.2]{Wbook}, we obtain the $\log$-Harnack inequality and Harnack inequality with power for $\bar{P}_{T}$. Next, we only show Theorem \ref{T2.2} (2),  and (1) is completely similar.

For every $p>(1+K_{2}K_{3})^{2}$, the Harnack inequality with power
\beq\label{5.41}
\bar{P}_{T}f(\eta)\leq \left(\bar{P}_{T} f^{p}(\xi)\right)^{\frac{1}{p}}\exp{\tilde{\Phi}_{p}(T;\xi,\eta)},\ \  \xi,\eta\in\C
\end{equation}
holds for non-negative function $f\in\B_{b}(\C)$, where
\begin{equation*}
\tilde{\Phi}_{p}(T;\xi,\eta)=\tilde{C}(p)\left\{1+\frac{|\xi(0)-\eta(0)|^2}{T-r}+\|\xi-\eta\|_{\infty}^{2}\right\}.
\end{equation*}
and $\tilde{C}:((1+K_{2}K_{3})^{2},\infty)\to(0,\infty)$ is a decreasing function.

From \eqref{5.12}, we obtain that for every $p>(1+K_{2}K_{3})^{2}$ and non-negative function $f\in\B_{b}(\C)$,
\beq\beg{split}\label{5.42}
P_{T}f(\eta)&=\bar{P}_{T} (f \circ\theta_{T}^{-1})(\theta_{0}(\eta))\\
&\leq\{\bar{P}_{T}(f^{p}\circ\theta_{T}^{-1})(\theta_{0}(\xi))\}^{\frac{1}{p}}\exp{\tilde{\Phi}_{p}(T;\theta_{0}(\xi),\theta_{0}(\eta))}\\
&=\{P_{T}f^{p}(\xi)\}^{\frac{1}{p}}\exp{\tilde{\Phi}_{p}(T;\theta_{0}(\xi),\theta_{0}(\eta))}\quad \xi,\eta\in\C
\end{split}\end{equation}
Combining \eqref{5.2}, we have
\begin{equation*}
\tilde{\Phi}_{p}(T;\theta_{0}(\xi),\theta_{0}(\eta))\leq \frac{81}{64}\tilde{C}(p)\left\{1+\frac{|\xi(0)-\eta(0)|^2}{T-r}+\|\xi-\eta\|_{\infty}^{2}\right\}=:\Psi_{p}(T;\xi,\eta).
\end{equation*}
Letting $K=K_1K_2$, $C=\frac{81}{64}\tilde{C}$, \eqref{5.42} yields \eqref{2.8}. Thus, we finish the proof.

\end{proof}

\paragraph{Acknowledgement} The authors would like to thank Professor Feng-Yu Wang and Jianhai Bao for corrections and helpful comments.

\beg{thebibliography}{99}

\bibitem{W} F.-Y. Wang,  \emph{Gradient estimate and applications for SDEs in Hilbert space with multiplicative noise and Dini continuous drift,}  J. Differ. Equations, 260 (2016), 2792-2829.

\bibitem{WH} F.-Y. Wang, X. Huang \emph{Functional SPDE with Multiplicative Noise and Dini Drift,}  to appear in Toulouse Sci. Math.

\bibitem{DZ} G. Da Prato, J. Zabczyk, \emph{Stochastic Equations in Infinite Dimensions,} Cambridge University Press, Cambridge, 1992.

\bibitem{G} T. E. Govindan, \emph{Mild Solutions of Neutral Stochastic Partial Functional Differential Equations,} International Journal of Stochastic Analysis. 2011, 186206.

\bibitem{LW08} W. Liu, F.-Y. Wang, \emph{Harnack inequality and strong Feller property for stochastic fast-diffusion equations,} J. Math. Anal. Appl. 342(2008), 651--662.

\bibitem{PR07} C. Pr\'ev\^ot, M. R\"ockner, \emph{A Concise Course on Stochastic Partial Differential Equations,} Lecture Notes in Math. vol. 1905, Springer 2007, Berlin.

\bibitem{PW} E. Priola, F.-Y. Wang, \emph{Gradient estimates for diffusion semigroups with singular coefficients,}  J. Funct. Anal. 236(2006), 244--264.

\bibitem{RW10} M. R\"ockner, F.-Y. Wang, \emph{Harnack and functional inequalities for generalized Mehler semigroups,}  J. Funct. Anal.  203(2007), 237--261.

\bibitem{W07} 	F.-Y. Wang, \emph{Harnack inequality and applications for stochastic generalized porous media equations,}  Annals of Probability 35(2007), 1333--1350.

\bibitem{Wbook} F.-Y. Wang, \emph{Harnack Inequality and Applications for Stochastic Partial Differential Equations,} Springer, New York, 2013.

\end{thebibliography}

\end{document}